\documentclass{amsart}

\usepackage{amssymb,amsmath,amsthm,delarray,stmaryrd}

\renewenvironment*{proof}{\noindent \textit{ Proof.
}}{ \makebox[11.9cm][r]{$\square$} }

\newtheorem*{theorem}{Theorem}
\newtheorem{lemma}{Lemma}
\newtheorem*{cor}{Corollary}

\newcommand*{\Y}{\mathbin{\mathsf Y}}

\renewenvironment*{proof}{\noindent \textit{ Proof.
}}{ \makebox[11.9cm][r]{$\square$} }

\newcommand*{\gp}[1]{\langle#1\rangle}

\makeatletter
\def\@setsubjclass{%
   2000 {\itshape Mathematics subject classification\/}%
   :\enspace\@subjclass\@addpunct.}
\def\@setkeywords{%
   {\itshape Keywords and phrases}:\enspace \@keywords\@addpunct.}
\makeatother

\begin{document}

\title[Group algebras with unit group of class $p$]
{Group algebras with unit group of class $p$}

\author{ Zs. Balogh}
\address{Institute of Mathematics and Informatics,
College of Ny\'\i regyh\'aza,  S\'ost\'oi \'ut 31/b, H-4410 Ny\'\i
regyh\'aza, Hungary} \email{baloghzs@nyf.hu}

\author{A. Bovdi}
\address{Institute of Mathematics, University of Debrecen, P.O. Box 12, 4010 Debrecen,
Hungary}
\email {bodibela@math.klte.hu}



\begin{abstract}
Let $V({\Bbb F}_{p}G)$ be the  group of normalized units of the
group algebra  ${\Bbb F}_{p}G$ of  a finite nonabelian $p$-group
$G$  over the field ${\Bbb F}_{p}$ of $p$ elements. Our goal is to
investigate the power structure of $V({\Bbb F}_{p}G)$, when it has
nilpotency class $p$. As a consequence,  we have proved that if
$G$ and $H$ are $p$-groups with cyclic Frattini subgroups and $p>
2$, then $V({\Bbb F}_{p}G)$ is isomorphic to $V({\Bbb F}_{p}H)$ if
and only if  $G$ and $H$ are isomorphic.
\end{abstract}

\thanks{
Supported by (OTKA) grants No.~T037202 and No.~T043034}

\maketitle

\section{Introduction}

Let $G$  be a finite  $p$-group  and ${\Bbb F}_{p}G$ its group
algebra over the field  ${\Bbb F}_{p}$ of $p$ elements. The
subgroup
$$
V({\Bbb F}_{p}G)=\Bigl\{ \sum_{g \in G }\alpha_g g \in {\Bbb
F}_{p}G\, \Bigm | \, \sum_{g \in G }\alpha_g=1 \Bigr\}
$$
is called {\it the group of normalized units}. Evidently $V({\Bbb
F}_{p}G)$ is a finite $p$-group and its order is $p^{|G|-1}$.

It is well-known that if $G$ is a finite nonabelian $p$-group then
the nilpotency class of  $ V({\Bbb F}_{p}G)$ is at least $p$. Mann
and Shalev \cite{mann} have recently shown that the nilpotency
class of $ V({\Bbb F}_{p}G)$ is  $p$ if and only if  the
commutator subgroup $G'$ is of order $p$. Our goal here is to
study the power structure of  $ V({\Bbb F}_{p}G)$ with nilpotency
class $p$ and to apply these results  to the isomorphism problem
of the group of units.

Since the modular isomorphism problem for finite $p$-groups has
been positively resolved  for some classes  of finite $p$-groups,
it is natural to investigate the question, originally raised by
S.~D.~Berman,  whether for a field  ${\Bbb F}_{p}$ of
characteristic $p$ and finite $p$-groups $G$ and $H$ the
isomorphism of  groups $V({\Bbb F}_{p}G)$ and $V({\Bbb F}_{p}H)$
implies the isomorphism of $G$ and $H$.

Note that Deskins \cite{deskins} gave a positive answer to the
isomorphism problem for finite abelian groups and  the Berman's
conjecture is true for these groups.

Further, the authors in \cite{balbod} solved this problem for the
group algebras of $2$-groups of maximal class over the field of
two elements, which is an extension of Baginski's result
\cite{baginski}.

Note that, in general, for nonabelian $p$-groups the Berman's
question is still  open. In the present paper we   give  a
positive answer to Berman's conjecture for $p$-groups of odd order
with cyclic Frattini subgroup.

\bigskip
We introduce some additional notation. Let $\zeta(G)$ be the
center and $\Phi(G)$ the Frattini subgroup of $G$,
 respectively;   $G^p=\gp{\;g^p\; | \; g \in G
\;}$, and let the commutator subgroup of $G$ be denoted by $G'$.
For each subset $D$ of $G$, let $\widehat{D}$ denote the element
$\sum_{g \in D} g$ of ${\Bbb F}_{p}G$.
\bigskip

\section{Results}

We begin with the description of the center $\zeta\big(V({\Bbb
F}_{p}G)\big)$ of $V({\Bbb F}_{p}G)$, where  $G$ is a finite
$p$-group with commutator subgroup $G'$ of order $p$ and $p>2$.
Let $C_{g_1},\ldots, C_{g_t}$ be all the different conjugacy
classes of $G$ which contain at least two elements. It is easy to
check that $t=\textstyle\frac{|G|-|\zeta(G)|}{p}$,
$$
\widehat{C_{g_i}}=g_i\widehat{G'}, \qquad
\widehat{C_{g_i}}\widehat{C_{g_j}}=0\qquad (1\leq i,j\leq  t)
$$
and   $\widehat{G'}$ is a central element with square $0$.

Clearly,  the set of all elements of the form $\sum_{i=1}^{t}
\alpha_{i} \widehat{C_{g_i}}$ is an ideal of the center
$\zeta({\Bbb F_p}G)$.  It follows that every central unit $x\in
\zeta\big(V({\Bbb F_p}G)\big)$ can be written as
\begin{equation}
x=z+\sum_{i=1}^{t} \alpha_{i}
\widehat{C_{g_i}}=z\Bigl(1+\sum_{i=1}^{t} \beta_{i}
\widehat{C_{g_i}}\Bigr)=z\prod_{i=1}^{t}\bigl(1+
g_i\widehat{G'}\bigr)^ { \beta_{i}},  \label{eq:egy}
\end{equation}
where $z\in V\big({\Bbb F_p}\zeta(G)\big)$ and $\alpha_{i},
\beta_{i}\in {\Bbb F_p}$.

Also it is easy to check that
$$
|V\big({\Bbb F_p}\zeta(G)\big)|=p^{|\zeta(G)|-1},
$$
\begin{equation} \big|\zeta\big(V({\Bbb F}_{p}G)\big)\big|=
p^{\frac{|G|+(p-1)|\zeta(G)|-p}{p}}, \label{eq:ketto}
\end{equation}
and according to (\ref{eq:egy}) we get
\begin{equation}
\zeta\big(V({\Bbb F_p}G)\big)=V\big({\Bbb F_p}\zeta(G)\big)\times
N, \label{eq:harom}
\end{equation}
where $ N=\prod_{i=1}^t \limits \gp{1+\widehat{C_{g_i}}}$ is an
elementary abelian subgroup  of $V({\Bbb F_p}G)$.

 Further, the commutator subgroup
$V({\Bbb F}_{p}G)'$ has exponent $p$ because it is a subgroup of
$1+{\frak I}(G')$ and ${\frak I}(G')^p=0$, where ${\frak I}(G')$
is the ideal generated by the elements of the form $h-1$ with
$h\in G'$.

Denote by $[{\Bbb F}_pG,{\Bbb F}_pG]$ the span of all ring
commutators $xy-yx$ with $x,y \in {\Bbb F}_pG$, which is called
\emph{ the commutator subspace} of ${\Bbb F}_pG$.

\begin{lemma} \label{lemma:submodule}\textup{(\textbf{Brauer},
\cite{brauer})} An element $\sum_{g \in G} \alpha_g g$ belongs to
the commutator subspace $[{\Bbb F}_pG,{\Bbb F}_pG]$ if and only if
$\sum_{g \in C_h} \alpha_g=0$ for every conjugacy class $C_h$ of
$G$. Moreover
$$
(x+y)^p \equiv x^p+y^p \pmod{[{\Bbb F}_pG,{\Bbb F}_pG]}
$$
for any $x,y \in {\Bbb F}_pG$.

\end{lemma}

We start by investigating the $p\,$th powers of the elements of $
V({\Bbb F}_{p}G)$.

\begin{lemma} \label{lemma:center} Let $G$ be a finite $p$-group with commutator
subgroup $G'$ of order $p>2$. Then $V({\Bbb F}_{p}G)^p$ is a
subgroup of the center $\zeta\big(V({\Bbb F}_{p}G)\big)$.
\end{lemma}

\begin{proof}
Let $H$ be the subgroup of $V({\Bbb F}_{p}G)$ generated by $x$ and
$y$, where $x\in V({\Bbb F}_{p}G)$, $g \in G$, and $y=g^{-1}xg$.
Evidently,
\[
\begin{split}
(x,y)&=x^{-1}(g,x)x(x,g)\\
     &=\big(x,(x,g)\big) \in
\gamma_3\big(V({\Bbb F}_{p}G)\big),
\end{split}
\]
so the subgroup $H\gamma_3\big(V({\Bbb
F}_{p}G)\big)/\gamma_3\big(V({\Bbb F}_{p}G)\big)$ is abelian. Thus
$H'$ is contained in $\gamma_3\big(V({\Bbb F}_{p}G)\big)$ and the
nilpotency class of $H$ is less than $p$. Then  $H$ is a regular
$p$-group and so, according to Theorem $12.4.2$ in \cite{hall}, we
have
$$
x^{-p}y^p=(x^{-1}y)^pd^p=(x,g)^pd^p
$$
for some element $d$ of the commutator subgroup of $V({\Bbb
F}_{p}G)$. But $V({\Bbb F}_{p}G)'$ has exponent $p$, so
$x^{-p}y^p=1$ and $x^p=g^{-1}x^pg$ for all $g \in G$. Thus $x^p$
is central, as asserted.

\end{proof}

\begin{lemma}\label{lemma:abp}
Let $H$ be a group generated by two elements $a,b$, and suppose
that its commutator subgroup $H'$ is central of prime order $p$.
In any group ring of $H$,
  \begin{equation}
    (a+b)^p=a^p+b^p+\sum_{r=1}^{p-1} \frac{1}{p}\binom{p}{r}
    a^rb^{p-r}\widehat{H'}. \label{eq:het}
  \end{equation}
\end{lemma}

\begin{proof}
As a first step, observe that $\zeta(H)=\gp{a^p,b^p}H'$ has index
$p^2$, and that the centralizer of any non-central element $h$ is
$\gp{h}\zeta(H)$. For $k,r \in \{ 1,2,\dots,p-1 \}$ and
$c_1,\dots, c_p\in \{a,b\}$ with $c_1\cdots c_p\in a^rb^{p-r}H'$,
it follows that no element of the coset $a^rb^{p-r}H'$ can commute
with the product $c_1\cdots c_k$ and so
\[
  c_{k+1}\cdots c_pc_1\cdots c_k=(c_1\cdots c_k)^{-1}(c_1\cdots c_p)(c_1\cdots
  c_k)\ne c_1\cdots c_p.
\]

Next, consider the set of all words $z_1z_2\cdots z_p$ of length
$p$ in the alphabet $\{ x,y \}$, as elements of the free semigroup
S freely generated by $\{ x,y \}$. The group of order $p$ acts on
this set by cyclically permuting the letters of a word. It is easy
to see that there are only two words fixed under this action,
$x^p$ and $y^p$. Since $p$ is prime, the length of each
non-singleton orbit is $p$. There are precisely $\binom{p}{r}$
words in which $x$ occurs $r$ times and $y$ occurs $p-r$ times,
and we conclude that these are permuted in
$\frac{1}{p}\binom{p}{r}$ orbits.

Let $\sigma : S \rightarrow H$ be the homomorphism defined by $x
\mapsto a,\; y \mapsto b$. The images of the orbits we counted all
lie in the coset $a^rb^{p-r}H'$. The point of the first step of
our argument was to show that the restriction of $\sigma$ to each
of these orbits is one-to-one. Since each orbit has length $p$ and
this is also the number of elements in the coset, it follows that
each element of $a^rb^{p-r}H'$ is the image of precisely
$\frac{1}{p}\binom{p}{r}$ of the words under consideration.

The rest of the proof may now be left to the reader.

\end{proof}

Let $x=\sum_{g \in G} \alpha_g g\in V({\Bbb F}_{p}G)$. From Lemma
\ref{lemma:submodule}, we know that $x^p=y+u$, where $y=\sum_{g
\in G}\alpha_g g^p$ and $u \in [{\Bbb F}_{p}G,{\Bbb F}_{p}G]$.
Lemma \ref{lemma:center} tells us that $x^p$ and the $g^p$ are
central in ${\Bbb F}_{p}G$; therefore so is $y$, and then also
$u$. By Lemma \ref{lemma:submodule} again, the support of $u$
cannot contain any element of $\zeta(G)$, so $u$ must be a linear
combination of the $\widehat{C_{g_i}}$. Thus $u^2=0$, and then
\begin{equation}
x^{p^2}=(y+u)^p=y^p=\Bigl( \sum_{g\in G} \alpha_g g^p
\Bigr)^p=\sum_{g \in G} \alpha_g g^{p^2}\label{eq:negy}
\end{equation}
(because each $g^p$ is central). We have proved that
$ V({\Bbb F}_{p}G)^{p^2}= V({\Bbb F}_{p}G^{p^2})$. When
$\exp{G}>p$, this shows that the exponents of the groups $V({\Bbb
F}_{p}G)$ and $G$ coincide. Consider next the case $\exp(G)=p$.
Choose $a,b\in G$ with $(a,b)\ne 1$. Then $a+b-1 \in V({\Bbb
F}_{p}G)$ and by (\ref{eq:het}) we have $(a+b-1)^p=(a+b)^p-1\ne
1$. It follows that
\begin{equation}
 \exp\big(V({\Bbb F}_{p}G)\big)=
 \begin{cases}
\exp(G)& \text{if} \quad \exp(G)>p;\\
p^2& \text{if}\quad \exp(G)=p.\label{eq:hat}
\end{cases}
\end{equation}
\begin{lemma}\label{lemma:huppert} \textup{(\textbf{Lemma III.9.6 in Huppert \cite{huppert1}})}
Let $U({\Bbb F}_p)$ be the group of units of ${\Bbb F}_p$ with odd
prime $p$. Then
\[
\sum_{\gamma\in U({\Bbb F_p})} \gamma^r=
\begin{cases} 0& \text{for}\quad 1\leq r\leq p-2;\\
                             p-1& \text{for}\quad r=p-1.
                             \end{cases}
\]
\end{lemma}

\begin{lemma}\label{lemma:vp} Let $G$ be a finite nonabelian $p$-group with
$|\Phi(G)|=p$.  Then $ V({\Bbb F}_{p}G)^p = V\big({\Bbb
F}_{p}G^p\big) \times N$, where $N=\prod_{i=1}^t \limits
\gp{1+\widehat{C_{g_i}}}$.
\end{lemma}

\begin{proof}
First we shall prove that $N\subseteq V({\Bbb F_p}G)^p$. Let
$\gamma\in U({\Bbb F}_p)$,\quad  $g\in G\setminus \zeta(G)$ and
$h\in G$ such that  $(g,h)\ne 1$. The commutator subgroup of
$\gp{h,g^{-1}h}$ coincides with $G'$ and
\begin{equation}
\widehat{G'}h^{-p}=\widehat{G'},\qquad  \big( (g^{-1}h)^p-1
\big)\widehat{G'}=0,\label{eq:nyolc}
\end{equation}
because $h^p,\, (g^{-1}h)^p \in G'$.

Clearly, for each $\gamma\in U({\Bbb F}_p)$ the element
$u_{\gamma}=h+\gamma (g^{-1}h-1)$ is a unit in ${\Bbb F}_pG$. By
(\ref{eq:het}) and (\ref{eq:nyolc})
\[
  \begin{split}
      u_{\gamma}^p&=\big( (h+\gamma g^{-1}h)-\gamma \big)^p=(h+\gamma g^{-1}h)^p-\gamma^p \\
                  &=h^p+\gamma\big( (g^{-1}h)^p -1 \big)+
                  \sum_{r=1}^{p-1}\textstyle \frac{1}{p}\binom{p}{r}h^r(\gamma
                  g^{-1}h)^{p-r}\widehat{G'}.
  \end{split}
\]

It follows that, with $\gamma$ ranging over $U({\Bbb F}_p)$,

\[
  \begin{split}
     \prod_{\gamma\in U({\Bbb F}_p)} (u_{\gamma}^p h^{-p})=1&+\Big( \sum_{\gamma\in U({\Bbb F}_p)}
     \gamma
       \Big)\big( (g^{-1}h)^p-1 \big)h^{-p} \\
       &+\sum_{r=1}^{p-1}\textstyle{\frac{1}{p}\binom{p}{r}}\Big( \sum\limits_{\gamma\in U({\Bbb F}_p)} \gamma^{p-r} \Big) h^r(
                  g^{-1}h)^{p-r}\widehat{G'}
  \end{split}
\]

and here, by Lemma \ref{lemma:huppert}, all summands with $r>1$
vanish, leaving
\[
 \prod_{\gamma\in U({\Bbb F}_p)} (u_{\gamma}^p h^{-p})=1-h(g^{-1}h)^{p-1}
 \widehat{G'}.
\]
Since $(g^{-1}h)^p\widehat{G'}=\widehat{G'}$ by (\ref{eq:nyolc}),
we have $h(g^{-1}h)^{-1}
 (g^{-1}h)^{p}\widehat{G'}=g\widehat{G'}$ and
\[
\prod_{\gamma\in U({\Bbb F}_p)} (u_{\gamma}^p
h^{-p})=1-g\widehat{G'}=(1+g\widehat{G'})^{-1}.
\]

Thus $1+g\widehat{G'}=\Big( \prod_{\gamma\in U({\Bbb F}_p)}
(u_{\gamma}^p h^{-p})\Big)^{-1} \in V({\Bbb F}_{p}G)^p$. Since the
elements of the form $1+g\widehat{G'}$ constitute a generator
system of $N$, we have proved that $N\subseteq V({\Bbb F_p}G)^p$,
as required.

Let $G^p$ be a nontrivial subgroup of $G$. Since $\Phi(G)$ is
cyclic, then $G^p=\gp{g^p}$ for some $g\in G$ and
$$
V({\Bbb F}_{p}\gp{g^p})\subseteq V({\Bbb F}_{p}\gp{g})^p\subseteq
V({\Bbb F}_{p}G)^p.
$$
Thus we have proved that $V\big({\Bbb F}_{p}G^p\big) \times N
\subseteq V({\Bbb F}_{p}G)^p$.

 Finally, the relation   $V({\Bbb F}_{p}G)^p \subseteq
V\big({\Bbb F}_{p}G^p\big) \times N$ follows from (\ref{eq:negy})
and the  prove is complete.

\end{proof}

The question: for which $p$-group $G$ is it true that $G\cap
V({\Bbb F}_{p}G)^{p}=G^{p}$, is due to Johnson \cite{jonhson}. The
previous lemma can be applied to conclude the following

\begin{cor} Let $G$ be a finite $p$-group such that  $|\Phi(G)|=p>2$. Then
$$
G\cap V({\Bbb F}_{p}G)^{p}=G^{p}.
$$
\end{cor}

\begin{proof} By Lemma \ref{lemma:vp} we get   $ V({\Bbb F}_{p}G)^p=V\big({\Bbb
F}_{p}G^{p}\big) \times N$, and  so
$$
G\cap V({\Bbb F}_{p}G)^{p}=G\cap V\big({\Bbb
F}_{p}G^{p}\big)=G^{p}.
$$

\end{proof}

\begin{theorem} Let $G$ and $H$ be finite nonabelian $p$-groups
with cyclic Frattini subgroup and $p>2$. Then $V({\Bbb F}_{p}G)$
is isomorphic to $V({\Bbb F}_{p}H)$ if and only if $G$ and $H$ are
isomorphic.
\end{theorem}

\begin{proof}
We say that $G$ is a central product $G_1{\Y}G_2$ of its subgroups
$G_1$ and $G_2$ if the elements of $G_1$ and $G_2$ commute and
together generate $G$, and $G_1 \cap G_2$ is the center of one of
the factors $G_1$, $G_2$.

It follows (for example) from Theorem $2$ in \cite{kovacs} that,
when $p>2$, every finite nonabelian $p$-group $G$ with cyclic
Frattini subgroup may be written as

\begin{equation}
G=E\times (K \Y L), \label{eq:centralis}
\end{equation}
where $E$ is elementary abelian, $K$ is either of order $p$ or an
extraspecial group of exponent $p$, and $L$ is either nontrivial
cyclic or nonabelian with a cyclic maximal subgroup, that is, $L$
is either $C_{p^n}=\gp{\,a\, \vert \, a^{p^n}=1\,}$  with  $n \geq
1$ or
$$
M_{p^n}=\gp{\, a,b \, \vert \,
a^{p^{n-1}}=b^p=1,\,(a,b)=a^{p^{n-2}}\,}\qquad   \text{with}\quad
n \geq 3.
$$

It is obvious that $L$ is cyclic if and only if
$\exp(G)=\exp\big(\zeta(G)\big)$; in this case, $|L|=\exp(G)$ and
$|K|=p \cdot |G:\zeta(G)|$, and otherwise $|L|=p \cdot \exp(G)$
and $|L| =p^{-1}\cdot |G:\zeta(G)|$. Consequently, the isomorphism
type of such a group is determined by the orders and exponents of
the group and its center.

The nontrivial part of the proof of the theorem is the claim that
these four invariants of $G$ are recognizable from the isomorphism
type of $V({\Bbb F}_{p}G)$.

First, $|G|$ is recognizable from $|V({\Bbb F}_{p}G)|$, and then
$|\zeta(G)|$ can be computed from $\bigl|\zeta\bigl(V({\Bbb
F}_{p}G)\bigr)\bigr|$ and (\ref{eq:ketto}). Using (\ref{eq:egy}),
it is not hard to see that
$$
\exp(\zeta(G))=\exp\Big(\zeta\big(V({\Bbb F}_{p}G)\big)\Big).
$$

It remains to show that the exponent of $G$ is also recognizable.
By (\ref{eq:hat}) if $\exp\big(V({\Bbb F}_{p}G)\big)>p^2$, then
$\exp\big(V({\Bbb F}_{p}G)\big)= \exp(G)$ and so we only have an
issue when $\exp\big(V({\Bbb F}_{p}G)\big)=p^2$.

In this outstanding case $\exp(G)\leq p^2$, so we obtain that
$|\Phi(G)|=p$. Indeed, it is clearly for the group $G$ with
$\exp(G)=p$. If $\exp(G)$ is equal to $p^2$ then by
(\ref{eq:centralis}) $\exp(G)=\exp(L)=p^2$. Clearly
$\Phi(G)=\Phi(L)$ and from the fact that $L$ is isomorphic to
either $C_{p^2}$ or $M_{p^3}$ follows that the subgroup $\Phi(G)$
has order $p$.

Then by Lemma \ref{lemma:vp} $ V({\Bbb F}_{p}G)^p = V\big({\Bbb
F}_{p}G^p\big) \times N$, where
$|N|=p^{\frac{|G|-|\zeta(G)|}{p}}$. Therefore we can determine
whether $|G^p|$ is either $1$ or $p$. We have proved that  the
exponent of $G$ is also recognizable, as required.

\end{proof}

\setlength{\parindent}{0pt}

\end{document}